%% file: main.tex
\documentclass[a4paper]{ifacconf}

\usepackage{graphicx}      
\usepackage{natbib}        

\usepackage{amsmath,amssymb,amsfonts}
\usepackage{xfrac}
\usepackage{xspace}
\usepackage{bbm}		
\usepackage{array}
\usepackage{algorithm}
\usepackage{algpseudocode}
\usepackage{multirow}
\usepackage{graphicx}
\usepackage{textcomp}
\usepackage{booktabs}
\usepackage{subfig}
\usepackage{footnote}
\usepackage{threeparttable}
\def\BibTeX{{\rm B\kern-.05em{\sc i\kern-.025em b}\kern-.08em
        T\kern-.1667em\lower.7ex\hbox{E}\kern-.125emX}}

\usepackage{balance}
\usepackage{mathtools}
\usepackage{newfloat}
\usepackage{caption}
\usepackage{enumitem}

\usepackage{alphalph}
\usepackage{etoolbox}

\patchcmd{\subequations}{\alph{equation}}{\alphalph{\value{equation}}}{}{}

\DeclarePairedDelimiter{\abs}{\lvert}{\rvert}

\newcommand{\trans}[1]{#1'}

\newtheorem{definition}{Definition}
\newtheorem{subdefinition}{Definition}[definition]

\usepackage{balance}

\usepackage{tikz}
\usepackage{pgf}
\usepackage{pgfplots}
\pgfplotsset{compat=1.13}
\usepgfplotslibrary{groupplots}

\begin{document}
\begin{frontmatter}

\title{Towards a Framework for Nonlinear Predictive Control using Derivative-Free Optimization\thanksref{footnoteinfo}}
\thanks[footnoteinfo]{The support of the EPSRC Centre for Doctoral Training in High Performance Embedded and Distributed Systems (HiPEDS, Grant Reference EP/L016796/1) and a Royal Society International Exchanges Grant (IES/R3/170011) is gratefully acknowledged.}

\author[ICL_EEE]{Ian McInerney} 
\author[ICL_EEE]{Lucian Nita}
\author[ICL_Aero]{Yuanbo Nie}
\author[Genoa]{Alberto Oliveri}
\author[ICL_EEE,ICL_Aero]{Eric C.\ Kerrigan}

\address[ICL_EEE]{Department of Electrical \& Electronic Engineering, Imperial College London, SW7~2AZ London, UK, email: \{i.mcinerney17,lucian.nita16,e.kerrigan,\}@imperial.ac.uk}
\address[ICL_Aero]{Department of Aeronautics, Imperial College London, SW7~2AZ London, UK}
\address[Genoa]{Department of Electrical, Electronic, Telecommunications Engineering and Naval Architecture, University of Genoa, via Opera Pia 11a, 16145, Genova, IT, email: alberto.oliveri@unige.it}

\begin{abstract}
The use of derivative-based solvers to compute solutions to optimal control problems with non-differentiable cost or dynamics often requires reformulations or relaxations that complicate the implementation or increase computational complexity.
 We present an initial framework for using the derivative-free Mesh Adaptive Direct Search (MADS) algorithm to solve Nonlinear Model Predictive Control problems with non-differentiable features without the need for reformulation.
 The MADS algorithm performs a structured search of the input space by simulating selected system trajectories and computing the subsequent cost value.
 We propose handling the path constraints and the Lagrange cost term by augmenting the system dynamics with additional states to compute the violation and cost value alongside the state trajectories, eliminating the need for reconstructing the state trajectories in a separate phase.
 We demonstrate the practicality of this framework by solving a robust rocket control problem, where the objective is to reach a target altitude as close as possible, given a system with uncertain parameters.
 This example uses a non-differentiable cost function and simulates two different system trajectories simultaneously, with each system having its own free final time.
\end{abstract}


\end{frontmatter}

\section{Introduction}

Model Predictive Control (MPC) has grown in popularity, due to its explicit handling of system constraints and the availability of efficient optimization solvers to compute the control solution.
 With the rise of data-based control and economic MPC, more engineers are wanting to control systems described by higher-fidelity models that capture nonlinearities or systems that only have blackbox/data-based models.
 These systems may pose a challenge for solvers that use first or second-order methods, since these methods require derivative information for the system dynamics, which may not be available or easily attainable. 
 
MPC can also be used with derivative-free optimization methods, where no knowledge of the derivatives for the optimization problem is required.
 To solve the Nonlinear MPC (NMPC) problem, these methods perform simulations of the dynamics to locate the future input trajectory that minimizes the cost function.
 These solvers can lead to embarrassingly parallel implementations, since all simulations in an iteration can be run in parallel.
 Simulation-based solvers have been readily used in Finite Control Set (FCS) algorithms in power electronics \citep{kouroModelPredictiveControl2015}.
 The derivative-free methods used to solve the FCS MPC problem are inefficient for long time horizons though, since they usually perform an exhaustive search of the possible future input sequences, which leads to a combinatorial explosion in the search space size as the horizon grows.
 Prior work has suggested several ways to work around this combinatorial explosion by randomly sampling a fixed number of points in the search space \citep{Joos2012}, or using methods such as pattern search \citep{gibsonDirectSearchApproach2015}, the Nelder-Mead simplex algorithm \citep{Sadrieh2011_GPU_NelderMead}, Trust Region methods \citep{daehlenNonlinearModelPredictive2014}, or Particle Swarm Optimization \citep{Xu2016}.
 
In this work, we build on the pattern search method from~\citet{gibsonDirectSearchApproach2015} and instead propose using the Mesh Adaptive Direct Search (MADS) algorithm from \citet{audetMeshAdaptiveDirect2006}.
 Using MADS, the number of points evaluated around the current iterate grows linearly with the input dimension, instead of combinatorially.
 We also propose an NMPC formulation for derivative-free optimization solvers that provides an easy way to handle the problem's path constraints and Lagrange cost term.
 
In Section~\ref{sec:nmpc} we introduce the NMPC optimization problem we use throughout the rest of the work, and in Section~\ref{sec:mads} we present an overview of the MADS algorithm.
 In Section~\ref{sec:formulation}, we show how to transform the NMPC problem into a form suitable for the MADS algorithm.
 We then present an example showing MADS solving an NMPC problem in Section~\ref{sec:numerical}, and conclude the paper with some future research directions in Section \ref{sec:future}.

\section{Nonlinear Model Predictive Control}
\label{sec:nmpc}

NMPC can be formulated as the optimization problem
\begin{subequations}
    \label{eq:nmpc}
    \begin{align}
    \underset{x,u,t_{f}}{\text{min}}\   & \Phi(x(t_{f}), u(t_{f}), t_{f}) + \int_{t_{0}}^{t_{f}} L(x(t), u(t), t) dt
    \label{eq:nmpc:cost}\\
    \text{s.t.\ }
    & f(x(t), \dot{x}(t), u(t), t) = 0, \quad \forall t \in [t_{0}, t_{f}] \label{eq:nmpc:dyn} \\
    & g(x(t), u(t), t) \leq 0, \label{eq:nmpc:pathconstraints} \\
    & h(x(0), u(0), t_{0}, x(t_{f}), u(t_{f}), t_{f}) = 0, \label{eq:nmpc:boundaryconds} 
    \end{align}
\end{subequations}
where $x:[t_{0}, t_{f}]\rightarrow \mathbb{R}^{n_{x}}$ and $u:[t_{0}, t_{f}]\rightarrow \mathbb{R}^{n_{u}}$ are the continuous-time state and input trajectories, respectively.
$f$ is the continuous-time nonlinear system dynamics, $g$ is the path constraints, and $h$ is the boundary conditions.
The cost functional~\eqref{eq:nmpc:cost} is composed of two terms: the Mayer cost $\Phi$, and the Lagrange cost $L$.

\section{Mesh Adaptive Direct Search}
\label{sec:mads}

This section provides a tutorial on the Mesh Adaptive Direct Search (MADS) algorithm. We combine the ideas from several works into a single statement of the algorithm, and provide a thorough discussion on the two techniques for implementing constraints in MADS.
 The notation used in this section has been slightly modified from the original MADS papers to make it consistent and to allow for a clearer description of the NMPC framework in Section~\ref{sec:formulation}.
 
MADS  is an extension of the Generalized Pattern Search (GPS) derivative-free method, and was first proposed in \citet{audetMeshAdaptiveDirect2006}.
 MADS has been extended to work with other types of variables (such as periodic, granular, integer or binary), model-based techniques, and multi-objective optimization \citep{audetDerivativeFreeBlackboxOptimization2017}.
 The optimization problem that we solve using MADS is
    \begin{align*}
    \underset{c, w}{\text{min}}\   & \mathcal{F}(c)
    \\
    \text{s.t.\ }
     (c, w) \in \Omega \coloneqq & \{ c \in \mathbb{R}^{m}, w \in \mathbb{R}^j:  \omega_i(c, w) \leq 0, \forall i \in K \},
    \end{align*}
where $\mathcal{F} : \mathbb{R}^{m} \rightarrow \mathbb{R}$ is an arbitrary function that computes the cost of the optimization problem.
 The vector $c \in \mathbb{R}^{m}$ contains the optimization variables that form the \textit{search space} 
 and $w \in \mathbb{R}^{j}$ is a vector of internal variables computed by $\mathcal{F}$ and used only in the constraints.
 The constraint set $\Omega$ is defined by the functions $\omega_i(\cdot)$ spanning the search space and the internal variables of the cost function $\mathcal{F}$, with $K$ the set of indices for the functions $\omega_{i}(\cdot)$ that define 
 $\Omega$.

\subsection{MADS Algorithm}

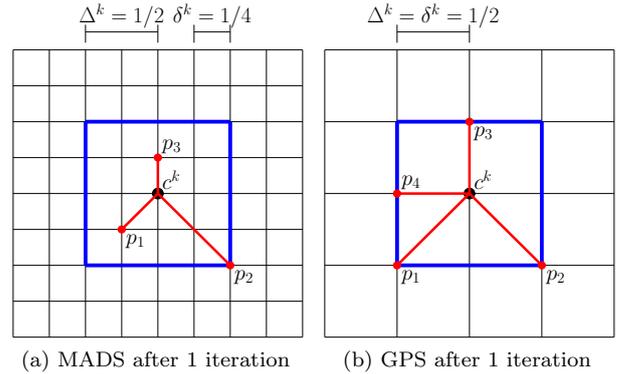
\begin{figure}[t]
    \centering
    \subfloat[MADS after 1 iteration\label{fig:mesh:mads}]{
    \resizebox{0.45\columnwidth}{!}{
    \input{figures/MADS_smallMesh.tex}
    }}
    \subfloat[GPS after 1 iteration\label{fig:mesh:gps}]{
    \resizebox{0.45\columnwidth}{!}{
    \input{figures/GPS_smallMesh.tex}
    }}

    \caption{The mesh (solid black lines spaced $\delta^k$ apart) and frame (region of size $\Delta^k$ inside the solid blue lines) for GPS and MADS after 1 successful iteration.}
    \label{fig:mesh}
\end{figure}

Pattern search methods sample the search space at a set of points called \textit{poll points} in each iteration.
 The poll points are located on a \textit{mesh} of size $\delta^{k}$ inside a \textit{frame} of size~$\Delta^{k}$ around the current iterate $c^{k}$, as in Fig.~\ref{fig:mesh}.
 The poll point with the lowest cost value is chosen as the next iterate.
 
The overall MADS algorithm is given in Algorithm~\ref{alg:mads}, and consists of three main parts: a search phase, a poll phase, and a mesh/frame adaptation phase.
 In the search phase, a set of points located on the mesh is generated (with no restriction on the method used to generate the points).
 The cost function is evaluated at these points, and if any have a lower cost value than the current iterate, the mesh and frame size are enlarged, and the search phase is repeated.
 If no better point is found, the poll phase is started.
 
In the poll phase, a set of polling directions that form a positive spanning set of the search space are generated to create a set of $m{+}1$ or more poll points contained inside the frame.
 The cost function is 
 evaluated at each poll point.
 MADS can perform these evaluations in parallel and in an opportunistic fashion, meaning it can evaluate multiple poll points at once and end the current iteration immediately after finding a point with lower cost value.

The final step in the MADS algorithm is to adapt the mesh and frame based on the results of the polling.
 If a lower cost has been found, then the mesh and frame get enlarged by the adjustment parameter $\tau$ to allow for a faster exploration of the search space.
 If no lower cost value is found, then the mesh and frame sizes are shrunk by~$\tau$ to concentrate the search closer to the current iterate.
 MADS can be terminated once the mesh size reaches a pre-determined threshold, meaning that there is no better point within that distance of the current iterate.

The poll and mesh adaptation phases are the only two phases required in every iteration of the MADS algorithm to get convergence.
 The search phase is optional; however, its use can speed up algorithm convergence, since it allows for more varied exploration of the search space.
 This can help the algorithm explore faster than the mesh can enlarge, and help it escape local minimums.

\subsection{Constraint Handling}

Implementing constraints inside  MADS  can be handled in at least one of two ways: 
    (1) directly constrain the poll points and either reject any that violate the constraints or modify the poll points, 
    (2) compute a metric indicating the amount of constraint violation and use it in the algorithm. 
If a constraint only contains variables in the algorithm's search space (e.g. only containing $c$ and not $w$), option 1 
allows for any poll points that violate the constraint to be skipped and not evaluated.
 Doing this can help to speed up the optimization process when $\mathcal{F}$ is expensive to evaluate, but does not provide any information to guide MADS from infeasible points to the feasible space.

To handle constraints that include both the poll point $c$ and internal variables $w$, the cost function can be redefined to be an extremal barrier function
\begin{equation}
    \label{eq:extremalBarrier}
    \mathcal{F}_{\Omega}(c) = \begin{cases}
    \mathcal{F}(c), & \forall c \in \Omega,\\
    \infty, &\forall c \notin \Omega.
    \end{cases}
\end{equation}
 The barrier term~\eqref{eq:extremalBarrier} ensures that MADS optimizes over only the feasible points by forcing all infeasible points to have infinite cost.
 While easy to implement, this
 has several drawbacks, including
    (i)~no information is provided to guide MADS from infeasible points to the feasible space,
    (ii)~a feasible starting point is required.

Alternately, \citet{audet2009progressive} proposed a progressive barrier technique that uses the constraint violation to keep both a feasible and infeasible iterate.
 In the poll phase, the $n_p$ polling directions are used to generate a set of $n_p$ poll points around both iterates ($2n_p$ total poll points).
 By tracking the best infeasible iterate, MADS can overcome the drawbacks of the extremal barrier method. 

We define the set $\Psi \supset \Omega$ as the relaxed constraint set, where some constraints are implemented using a progressive barrier form and others use an extremal barrier form.
 The function $\mathcal{H} : \mathbb{R}^{k+j} \rightarrow \mathbb{R}$ is then used to compute the constraint violation, and is defined as
 \begin{equation*}
    \mathcal{H}(c, w) =
    \begin{cases}
        0, & \text{if } (c,w) \in \Omega, \\
        \sum_{l \in K_{\Psi}} \left(\max\{\omega_{l}(c, w), 0\}\right)^2, & \text{if } (c,w) \in \Psi \setminus \Omega, \\
        \infty, & \text{otherwise},
    \end{cases}
\end{equation*}
where $K_{\Psi}$ is the set of indices for the constraints $\omega(\cdot)$ that use the progressive barrier form.

Inside MADS, the progressive barrier relies on the notion of a dominating iterate, which is a point that is better than the others in a specific sense.
\stepcounter{definition}
\begin{subdefinition}[Dominating Feasible Point]
    \label{def:point:feasDominating}
    The feasible point $p \in \Omega$ \textit{dominates} the point $y \in \Omega$, denoted $p \prec_{\mathcal{F}} y$, when $\mathcal{F}(p) < \mathcal{F}(y)$.
\end{subdefinition}
\begin{subdefinition}[Dominating Infeasible Point]
    \label{def:point:infeasDominating}
    The infeasible point $p \in \Psi \setminus \Omega$ \textit{dominates} the point $y \in \Psi \setminus \Omega$, denoted $p \prec_{\mathcal{H}} y$, when $\mathcal{F}(p) < \mathcal{F}(y)$ and $\mathcal{H}(p, w_{p}) \leq \mathcal{H}(y, w_{y})$ (with at least one strict inequality).
\end{subdefinition}
For feasible points, dominating means that the point has the smallest cost value, while for infeasible points it means that the point has both the smallest cost value and the smallest constraint violation.

An iteration of the MADS algorithm is one of three types: \textit{dominating}, \textit{improving} or \textit{unsuccessful}.
 A dominating iteration is when MADS finds a new point that has both a lower cost value and also a lower (possibly 0) constraint violation, so the iterates are updated, and the mesh is expanded to search in a larger region around the new iterates.
 In an improving iteration, no points with a lower cost were found, but an infeasible point with a constraint violation smaller than the penalty parameter $\eta$ was found.
 In this case the new point becomes the next infeasible iterate, the mesh  remains unchanged, and $\eta$ is set to the constraint violation at the new infeasible iterate.
 All other iterations are considered unsuccessful, so the mesh will be shrunken to search in a closer region to the iterates in the next iteration.
 Formal definitions for these iteration types can be found in \citet{audet2009progressive}.

The penalty parameter $\eta$ is nonincreasing with  iteration number, starting at $\infty$ and decreasing to $0$ as the algorithm runs.
 This means that at the beginning, MADS is prioritizing the search for the lowest cost value by allowing large constraint violations, but over time $\eta$ decreases and forces the infeasible iterates to move towards the feasible space.
 The penalty barrier constraint method behaves similar to a filter method in other optimization algorithms, such as SQP, but requires the barrier/penalty parameter to be only nonincreasing and not strictly decreasing.

\begin{algorithm}[pt!]
\begin{minipage}{\columnwidth}
\caption{Mesh Adaptive Direct Search with Progressive Barrier Constraints \citep{audetDerivativeFreeBlackboxOptimization2017}}
\label{alg:mads}
\begin{algorithmic}[1]
    \algnewcommand\And{\textbf{and}}
    \algnewcommand\Or{\textbf{or}\xspace}
    \algnewcommand\True{\textbf{True}}
    \algnewcommand\False{\textbf{False}}
    \algnewcommand\Goto{\textbf{goto}\xspace}
    \algnewcommand\Continue{\State \textbf{continue}}
    \algnewcommand\algorithmicinput{\textbf{Let:}}
    \algnewcommand\Let{\item[\algorithmicinput]}
    \algnewcommand\Blank{\Statex \vspace{-0.5em}}
    
    \Let $\mathbb{M}^{k}$ be the set of all mesh points using mesh size $\delta^{k}$
    \Let $\prec_{\mathcal{F}}$ and $\prec_{\mathcal{H}}$ be as given in Definitions~\ref{def:point:feasDominating} and~\ref{def:point:infeasDominating}
    \Require $\Delta^{0} \in (0, \infty)$ \Comment{Initial frame size}
    \Require $\tau \in (0, 1)$ \Comment{Mesh size adjustment parameter}
    \Require $\epsilon_{stop} \in [0, \infty)$ \Comment{Stopping tolerance}
    \Require $c_{f}^{0}$ \textbf{And\textbackslash Or} $c_{i}^{0}$  \Comment{Initial mesh centers}
    
    \Blank
    \Statex $k \gets 0$
    \While{ $\Delta^{k} \geq \epsilon_{stop}$ }\footnotemark[1]
        \State $\mathbb{V}^{k} \gets \emptyset$
        \State $\delta^{k} \gets \min{\{\Delta^{k}, (\Delta^{k})^2\}}$
        \Comment{Compute mesh size}
        
        \Blank
        \State \textbf{1) Search Phase:}
        \label{alg:mads:search}
        \State Generate search points $\mathbb{S}^{k} \subset \mathbb{M}^{k}$
        \ForAll{$s \in \mathbb{S}^{k}$}
        \State Compute $\mathcal{F}$ and $\mathcal{H}$ at $s$
        \If{ $s \prec_{\mathcal{F}} c_{f}^{k}$ }
            \State $c_{f}^{k+1} \gets s$
            \Comment{Dominating}
            \State \Goto line \ref{alg:mads:meshadaption}
        \ElsIf{ $s \prec_{\mathcal{H}} c_{i}^{k}$ }
            \State $c_{i}^{k+1} \gets s$, $\eta^{k+1} \gets \mathcal{H}(s)$
            \Comment{Dominating}
            \State \Goto line \ref{alg:mads:meshadaption}
        \EndIf
        \EndFor
        
        \Blank
        \State \textbf{2) Poll Phase:}
        \label{alg:mads:poll}
        \State Generate positive spanning set $\mathbb{D}^{k}$
        \State $\mathbb{P}^{k} \gets \{c_{f}^{k} + \delta^{k}d \quad \forall d \in \mathbb{D}^k\} \cup \{c_{i}^{k} + \delta^{k}d \quad \forall d \in \mathbb{D}^k\}$ 
        \ForAll{$p \in \mathbb{P}^{k}$}
        \State Compute $\mathcal{F}$ and $\mathcal{H}$  at $p$
        \If{ $p \prec_{\mathcal{F}} c_{f}^{k}$ }
            \State $c_{f}^{k+1} \gets p$
            \Comment{Dominating}
            \State \Goto line \ref{alg:mads:meshadaption}
        \ElsIf{ $p \prec_{\mathcal{H}} c_{i}^{k}$ }
            \State $c_{i}^{k+1} \gets p$, $\eta^{k+1} \gets \mathcal{H}(p)$
            \Comment{Dominating}
            \State \Goto line \ref{alg:mads:meshadaption}
        \ElsIf{ $\mathcal{H}(p) < \eta^{k}$ }
            \State $\mathbb{V}^{k} \gets \mathbb{V}^{k} \cup \{ p \}$
            \Comment{Improving}
        \EndIf
        \EndFor
        
        \Blank
        \If{ Improving }
            \State $v \gets \text{argmax}\{ \mathcal{H}(v) : \mathcal{H}(v) < \eta^{k}, v \in \mathbb{V}^{k}\}$
            \State $c_{i}^{k+1} \gets v$, $\eta^{k+1} \gets \mathcal{H}(v)$
        \EndIf
        
        \Blank
        \State \textbf{3) Mesh/Frame Update:}
        \label{alg:mads:meshadaption}
        \If{ Dominating }
            \State $\Delta^{k+1} \gets \tau^{-1} \Delta^{k}$
            \Comment{Expand frame size}
        \ElsIf{ Improving }
            \State $\Delta^{k+1} \gets \Delta^{k}$
            \Comment{Frame doesn't change}
        \Else
            \State $\Delta^{k+1} \gets \tau \Delta^{k}$
            \Comment{Shrink frame size}
        \EndIf
        
        \Blank
        \State $k \gets k+1$
    \EndWhile
\end{algorithmic}
\renewcommand{\thefootnote}{\arabic{footnote}}
\footnotetext{\footnotemark[1]If a variable is not set in an iteration, it retains the same value in the next iteration.}
\end{minipage}
\end{algorithm}

\subsection{Meshing}

In MADS, the frame and mesh are controlled by different parameters ($\Delta^k$ and $\delta^k$, respectively), with
\begin{equation}
    \delta^k = \min\{\Delta^k, (\Delta^k)^2\}
    \label{eq:madsUpdate}
\end{equation}
normally used. The main difference between MADS and GPS is that in GPS $\delta_k = \Delta_k$ at all times.
By allowing the mesh to shrink faster than the frame size, more 
polling points are created around the 
iterate $c^{k}$.
For example, Fig.~\ref{fig:mesh:mads} shows the result of one iteration of MADS using update rule \eqref{eq:madsUpdate}.
In this case there are 24 possible polling points for MADS, versus 8 for GPS in Fig.~\ref{fig:mesh:gps}.

\subsection{Convergence}

The fact that the mesh in MADS shrinks faster than the frame size means that, as the algorithm converges, there will be an asymptotically dense set of poll directions created.
It was shown in \citet{audetMeshAdaptiveDirect2006} that this set of asymptotically dense poll directions allows for the algorithm to converge to a Clarke stationary point with non-negative Clarke derivatives (e.g. a local minimum for the non-smooth function) when linear constraints are placed on the search space.
In contrast, GPS loses the theoretical convergence guarantees when simple bound constraints are applied.
Additionally, when a progressive barrier approach is used for handling nonlinear constraints, it was shown in \citet{audet2009progressive} that MADS is still effective at finding the stationary point.

\section{NMPC Problem Formulation for derivative-free Optimization}
\label{sec:formulation}

The core of the problem formulation we propose is a continuous-time shooting method based on \citet{gibsonDirectSearchApproach2015}, where the search space for MADS is the set of all input trajectories.
 Each function evaluation in the MADS algorithm simulates the system over the time horizon with the chosen input trajectory, and then computes the overall cost value and constraint violation for that trajectory.
 To easily handle the path constraints and the Lagrange term, we introduce an augmentation scheme where the system dynamics are augmented with new states representing the Lagrange term and the path constraint violation.
 
\subsection{Cost Functional}
\label{sec:formulation:cost}

 We split the NMPC cost functional~\eqref{eq:nmpc:cost} into its two terms, and augment the system dynamics with the Lagrange term by adding a new state $l(t)$ that represents the value of the Lagrange term at time $t$.
 This new state is governed by 
 \begin{equation}
    \label{eq:lagrangeDynamics}
     \dot{l}(t) = L(x(t), u(t), t),
 \end{equation}
 which is computed alongside the system's trajectory.
 
\subsection{Path Constraints}
\label{sec:formulation:pathConstraints}

The constraints~\eqref{eq:nmpc:pathconstraints} are formulated as progressive barrier constraints inside the set $\Psi$.
 An $L_1$ measure of constraint violation is used for each  constraint, meaning the value reported as the violation experienced by constraint $i$ is
 \begin{equation}
    \label{eq:pathConIntegral}
     v_{i} = \int_{0}^{t_{f}} \max\{ 0, g_i(x(t), u(t), t) \} dt.
 \end{equation}

To compute the integral~\eqref{eq:pathConIntegral}, we add new states $v(t)$ that represent the constraint violation over time, with those states being governed by the dynamics equation
\begin{equation}
    \label{eq:pathConDynamics}
    \dot{v}(t) = g^{+}(x(t), u(t), t)\\
\end{equation}
 where $g^{+}$ is the vector function representing the element-wise computation of $\max\{ 0, g_i(x(t), u(t), t) \}$.

\subsection{Overall Problem Formulation}
\label{sec:formulation:overall}

An implementation of the formulation described in this section is given in Algorithm~\ref{alg:nmpcFunc}.
 The first step is to construct the input trajectory described by the current poll point $c$.
 The structure of the input trajectory will depend on the problem, but possible options are:
     (i)~a zero-order hold trajectory with the input values at each sample given by the elements of the poll point,
     (ii)~an interpolated polynomial with the interpolation points 
     given by the elements of the poll point,
     (iii)~a feedback policy with the parameters of the policy given by the elements of the poll point.
 After constructing the input trajectory, the augmented system~\eqref{eq:mads:prob:diffeq} is simulated over the horizon length using a suitable numerical solver for differential equations.

After the simulation completes, the boundary condition violation is computed, and the total constraint violation is then computed using~\eqref{eq:mads:prob:progressiveBarrier}, where $\rho$ are weights that can be applied to various path and boundary constraints to give them more influence in the algorithm.
 Finally, the value of the NMPC cost function~\eqref{eq:nmpc:cost} for the selected poll point is computed using~\eqref{eq:mads:prob:cost}, which computes the Mayer term and then adds the final value of the Lagrange state $l(t_{f})$.
 

\begin{algorithm}[t]
\begin{minipage}{\columnwidth}
\caption{Function evaluation for the NMPC problem}
\label{alg:nmpcFunc}
\begin{algorithmic}[1]
 
    \algnewcommand\And{\textbf{and}}
    \algnewcommand\Or{\textbf{or}\xspace}
    \algnewcommand\True{\textbf{True}}
    \algnewcommand\False{\textbf{False}}
    \algnewcommand\Goto{\textbf{goto}\xspace}
    \algnewcommand\Continue{\State \textbf{continue}}
    \algnewcommand\algorithmicinput{\textbf{Let:}}
    \algnewcommand\Let{\item[\algorithmicinput]}
    \algnewcommand\Blank{\Statex \vspace{-0.5em}}
    
    \Let $c$ be the point in the search space being evaluated
    
    \State Construct the input trajectory $u$ from $c$
    \State Simulate the augmented dynamics~\eqref{eq:mads:prob:diffeq} using an appropriate solver for the differential equations
    \begin{align}
        \label{eq:mads:prob:diffeq}
        \begin{bmatrix}
            0 \\
            0 \\
            0
        \end{bmatrix} &=
        \begin{bmatrix}
            f(x(t), \dot{x}(t), u(t), t) \\
            L(x(t), u(t), t) - \dot{l}(t)\\
            g^{+}(x(t), u(t), t) - \dot{v}(t)
        \end{bmatrix}
    \end{align}
    
    \State Compute the violation of the boundary conditions
    \begin{equation}
        \label{eq:mads:prob:boundConstViol}
        v_{b} = \sum_{i} \rho_{b_i} \abs{h_{i}(x(0), u(0), t_{0}, x(t_{f}), u(t_{f}), t_{f})}
    \end{equation}
 
    \State Compute the overall constraint violation
    \begin{equation}
        \label{eq:mads:prob:progressiveBarrier}
        \mathcal{H} \gets v_{b}^{2} + \sum_{i} \rho_{i} ( v_{i}(t_{f}) )^2
    \end{equation}
    
    \State Compute the cost function value
    \begin{equation}
        \label{eq:mads:prob:cost}
        \mathcal{F} \gets \Phi(x(t_{f}), u(t_{f}), t_{f}) + l(t_{f})
    \end{equation}
\end{algorithmic}
\end{minipage}
\end{algorithm}

\subsection{Discussion of the Augmentation Scheme}

The proposed augmentation scheme provides an easy way to include constraints in the derivative-free optimization problem, but it may not be the most efficient way.
 Moving the constraints and the Lagrange term into the dynamics allows them to be on the same mesh as the dynamics, and removes the need for additional quadrature schemes to compute the final cost and constraint violation after the system has been simulated, as is usually done in other solvers.
 While this removes the need for an algorithm to handle the mesh refinement, it introduces a requirement that the dynamics solver must have an error control mechanism to limit the numerical integration error.
 
While easy to use, the augmentation scheme can also introduce inefficiencies to the solver.
 By introducing the dynamics equations~\eqref{eq:lagrangeDynamics} and~\eqref{eq:pathConDynamics}, we have made the solvability and the stiffness of~\eqref{eq:mads:prob:diffeq} be dependent on the Lagrange term and the constraints as well.
 This means that when the cost or constraints vary over time or are defined by stiff equations, they can cause the dynamics solver to struggle and possibly require smaller step sizes.

\section{Numerical Examples}
\label{sec:numerical}

We demonstrate the use of the formulation from Section~\ref{sec:formulation} by solving a robust rocket throttle control problem.
 We examine a laboratory scale rocket aiming to reach its apogee at a target altitude of 10,000 feet, governed by 
 \begin{equation*}
    \dot{x}(t) =
    f(x(t)) = 
    \begin{bmatrix}
       v_{v}(t)\\
       \frac{T(t)-0.5C_{D} \rho(h(t))\frac{\pi d^{2}}{4}({v_{v}(t)^{2}})}{m(t)}-g(h(t))\\
       -\frac{T(t)}{I_{sp}} \\
   \end{bmatrix},
 \end{equation*}
 with the state vector
 $x=\trans{\begin{bmatrix}
 h & v_{v} & m
 \end{bmatrix}}$
 composed of the rocket's altitude, vertical velocity and mass, respectively.
 $g(h)$ and $\rho(h)$ are the gravitational constant and the air density at altitude $h$, respectively, and $d$ is the rocket diameter.
 The model has an uncertain drag coefficient $C_D \in [ C_{D_{l}}, C_{D_{u}} ]$ and specific impulse $I_{sp} \in [ I_{sp_l}, I_{sp_u} ]$.
 The thrust is a function of time, $T(\cdot)$, that takes values from the set $T(t) \in [0, T_{max}]$.
 
It can be shown that the open-loop system is monotone with respect to the throttle setting and physical parameters ($C_D$, $I_{sp}$) \citep{angeliMonotoneControlSystems2003}, so all possible trajectories will lie in a tube defined by an upper $x_{u}$ and lower $x_{l}$ trajectory (using the upper and lower drag coefficients and specific impulses, respectively).
 We use a piecewise constant throttle function, and optimize over both the throttle settings $T_i$ and the switching times $\sigma_i$.
 A robust version of the problem then optimizes over the two trajectories $x_{l}$ and $x_{u}$, with the objective to minimize the largest deviation of the two trajectories from a target apogee in the min-max formulation 
 
 \newpage
 \begin{subequations}
 \label{eq:example:firstProblem}
 \begin{equation}
    \min_{T, \sigma}\   \max_{x_l,x_u} \left\{\abs{h_u(t_{f_u})-3048},\abs{3048-h_l(t_{f_l})}\right\} \label{eq:example:first:cost}
    \end{equation}
    \begin{align}
    \text{s.t.\ } & \dot{x}_l(t) = f(x_l(t)),\  
                   \dot{x}_u(t) = f(x_u(t)) 
                   \\
                  & x_l(0) = x_u(0) = \trans{\begin{bmatrix}0 & 0 & 33.5 \end{bmatrix}} \\
                  & v_{v_l}(t_{f_l}) = v_{v_u}(t_{f_u})=0 \\
                  & v_{v_l}(t) \leq 150 \qquad \forall t \in [0, t_{f_l}] \label{eq:example:first:lowerVelCon} \\
                  & v_{v_u}(t) \leq 150 \qquad \forall t \in [0, t_{f_u}] \label{eq:example:first:upperVelCon} \\
                  & m_l(t_{f_l}) \geq 26,\  m_u(t_{f_u}) \geq 26 \\
                  & 0 \leq \sigma_1 \leq \sigma_2 \leq \cdots \leq \sigma_i \\
                  & \sigma_i \leq \min{\{t_{f_l}, t_{f_u}\}} \\
                  & T_i \in [0, T_{max}] \qquad \forall i 
 \end{align}
 \end{subequations}
 where the same thrust function is used for both the upper and lower trajectories.
 The trajectories share the same thrust profile $T_{i}$ and switching times $\sigma_i$, but have their own final times given by $t_{f_u}$ and $t_{f_l}$ for the upper and lower trajectories respectively.
 As an example path constraint, we constrain the rocket's velocity to be less than 150\,m/s, representing a modeling constraint where the drag coefficient uncertainity bounds are only known below 150\,m/s.

\begin{figure}[t!]
  \centering
  \subfloat[Altitude profile\label{fig:goddard:constrained:altitude}]{
    \includegraphics[width=0.45\textwidth]{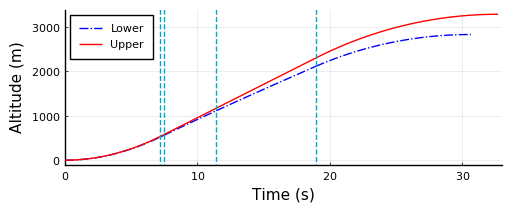}}\\
  \subfloat[Velocity profile\label{fig:goddard:constrained:velocity}]{
    \includegraphics[width=0.45\textwidth]{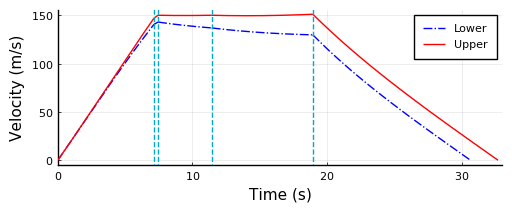}}\\
  \subfloat[Input trajectory\label{fig:goddard:constrained:input}]{
    \includegraphics[width=0.45\textwidth]{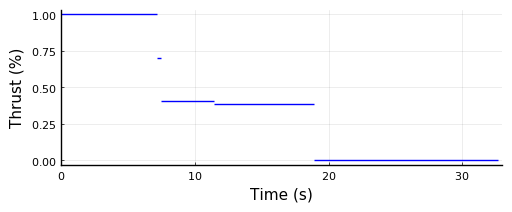}}
  \caption{Open-loop solution to~\eqref{eq:example:firstProblem} with 5 phases (Switching times are at the dotted lines).}
  \label{fig:goddard:constrained}
\end{figure}

MADS was used to  solve~\eqref{eq:example:firstProblem}  with a thrust trajectory composed of five intervals.
 This resulted in  optimal throttle settings shown in Fig.~\ref{fig:goddard:constrained:input}, which resemble the theoretically optimal `bang-singular-bang' solution.
 The resulting altitude and velocity profiles can be seen in Figures~\ref{fig:goddard:constrained:altitude} and~\ref{fig:goddard:constrained:velocity}, respectively.
 In the velocity profile, the constraint of 150\,m/s is satisfied by both the upper and lower trajectories when the constraint is implemented using additional states, as described in Section~\ref{sec:formulation:pathConstraints}.
 
\begin{figure}[t!]
  \centering
  \subfloat[No constraints enforced\label{fig:goddard:violation:unconstrained}]{
    \includegraphics[width=0.45\textwidth]{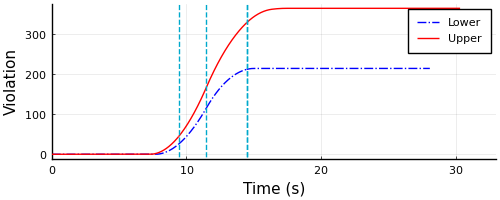}}\\
  \subfloat[Constraints enforced as progressive barrier constraints\label{fig:goddard:violation:constrained}]{
    \includegraphics[width=0.45\textwidth]{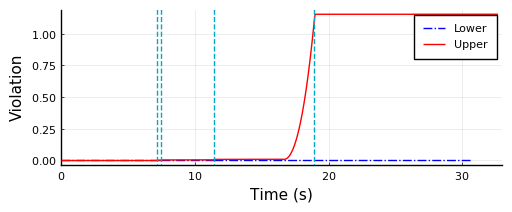}}
  \caption{Cumulative violation of constraints~\eqref{eq:example:first:lowerVelCon} and~\eqref{eq:example:first:upperVelCon}}
  \label{fig:goddard:violation}
\end{figure}

The evolution of the violation state for constraints~\eqref{eq:example:first:lowerVelCon} and~\eqref{eq:example:first:upperVelCon} can be seen in Fig.~\ref{fig:goddard:violation}.
 The violations experienced when the velocity constraints are not enforced 
 are shown in Fig.~\ref{fig:goddard:violation:unconstrained}, and are monotonically increasing along the horizon.
 Once the constraints are enforced, the constraint violation decreases significantly, as shown in Fig.~\ref{fig:goddard:violation:constrained}.
 
The non-zero violation for the upper trajectory 
is due to the current implementation of the progressive barrier constraints in the DirectSearch.jl\footnote{https://github.com/ImperialCollegeLondon/DirectSearch.jl} Julia package that was used.
 In the current implementation, there is no method for passing the internal variables $w$ from the computation of $\mathcal{F}$ to the computation of $\mathcal{H}$.
 This means that~\eqref{eq:mads:prob:progressiveBarrier} could not be used to compute the progressive barrier penalty value, so instead the velocity constraints were implemented by adding a penalty term in the cost for the deviation from zero of the final value of the velocity constraint states $v(t)$.
 
\section{Conclusions and Future Directions}
\label{sec:future}

We presented a formulation for NMPC that includes path constraints and the Lagrange cost term in a form suitable for use with MADS. 
 We have also shown that this can be used to solve problems with non-differentiable cost functions, discontinuous dynamics and free end times.
 
This formulation provides the initial steps to implement NMPC with MADS, but there are more open questions remaining.
 In the robust rocket control example in Section~\ref{sec:numerical}, we showed a zero-order hold input sequence, but this sequence can become impractical to optimize over for systems with a long horizon or many inputs.
 Additional work should be done to examine how to effectively optimize over the other input representations given in Section~\ref{sec:formulation:overall}, and to determine if there is a preferred representation to use.

The reliance on a single shooting method for simulating the system trajectories in Algorithm~\ref{alg:nmpcFunc} is also problematic.
 Small changes to the input sequence/initial conditions can drastically affect the simulated trajectory and the resulting constraint violation measure, so searching for an input trajectory that doesn't violate the boundary conditions or path constraints becomes more difficult.
 Future work should explore if a multiple shooting paradigm can be adapted into the proposed NMPC framework for MADS, with a particular focus on efficient handling of the defect constraints between the shooting regions.
 Additionally, switching to a multiple shooting approach may add opportunities for parallelism in the formulation itself by simulating all shooting intervals in parallel.

Finally, the example demonstrates this framework in use as an open-loop solver for optimal control problems, but the parallelization potential and simplicity of the computations may make the MADS-based framework beneficial for closed-loop control.
 More work should be done to examine how this framework behaves in a closed-loop implementation, and how a warm-started MADS implementation could be used inside a real-time iteration type framework.
 
\bibliography{references}

\end{document}

%% file: figures/MADS_smallMesh.tex
\begin{tikzpicture}

\tikzstyle{every node}=[font=\LARGE]

\node[draw, circle, inner sep = 3pt, fill] at (0, 0) {};
\node[above right] at (0, 0) {$c^{k}$};

\foreach \x in {-4, ..., 4}
{
	\draw (\x, -4) -- (\x, 4) {};
}

\foreach \y in {-4, ..., 4}
{
	\draw (-4, \y) -- (4, \y) {};
}

\draw[line width=3pt, blue] (2, 2)   -- (2, -2) {};
\draw[line width=3pt, blue] (2, 2)   -- (-2, 2) {};
\draw[line width=3pt, blue] (-2, -2) -- (2, -2) {};
\draw[line width=3pt, blue] (-2, -2) -- (-2, 2) {};

\node[draw, circle, inner sep = 2pt, fill, red] at (-1, -1) {};
\node[below right] at (-1, -1) {$p_{1}$};
\draw[line width = 2pt, red] (0, 0) -- (-1, -1) {};

\node[draw, circle, inner sep = 2pt, fill, red] at (2, -2) {};
\node[below right] at (2, -2) {$p_{2}$};
\draw[line width = 2pt, red] (0, 0) -- (2, -2) {};

\node[draw, circle, inner sep = 2pt, fill, red] at (0, 1) {};
\node[above right] at (0, 1) {$p_{3}$};
\draw[line width = 2pt, red] (0, 0) -- (0, 1) {};

\draw (-2, 4.5) -- (0, 4.5) {};
\draw (-2, 4.2) -- (-2, 4.7) {};
\draw (0, 4.2) -- (0, 4.7) {};
\node[above] at (-1, 4.5) {$\Delta^{k} = 1/2$};

\draw (1, 4.5) -- (2, 4.5) {};
\draw (2, 4.2) -- (2, 4.7) {};
\draw (1, 4.2) -- (1, 4.7) {};
\node[above] at (1.5, 4.5) {$\delta^{k} = 1/4$};

\end{tikzpicture}

%% file: figures/GPS_smallMesh.tex
\begin{tikzpicture}

\tikzstyle{every node}=[font=\Large]

\node[draw, circle, inner sep = 3pt, fill] at (0, 0) {};
\node[above right] at (0, 0) {$c^{k}$};

\foreach \x in {-4, -2, 0, 2, 4}
{
	\draw (\x, -4) -- (\x, 4) {};
}

\foreach \y in {-4, -2, 0, 2, 4}
{
	\draw (-4, \y) -- (4, \y) {};
}

\draw[line width=3pt, blue] (2, 2)   -- (2, -2) {};
\draw[line width=3pt, blue] (2, 2)   -- (-2, 2) {};
\draw[line width=3pt, blue] (-2, -2) -- (2, -2) {};
\draw[line width=3pt, blue] (-2, -2) -- (-2, 2) {};

\node[draw, circle, inner sep = 2pt, fill, red] at (-2, -2) {};
\node[below right] at (-2, -2) {$p_{1}$};
\draw[line width = 2pt, red] (0, 0) -- (-2, -2) {};

\node[draw, circle, inner sep = 2pt, fill, red] at (2, -2) {};
\node[below right] at (2, -2) {$p_{2}$};
\draw[line width = 2pt, red] (0, 0) -- (2, -2) {};

\node[draw, circle, inner sep = 2pt, fill, red] at (0, 2) {};
\node[below right] at (0, 2) {$p_{3}$};
\draw[line width = 2pt, red] (0, 0) -- (0, 2) {};

\node[draw, circle, inner sep = 2pt, fill, red] at (-2, 0) {};
\node[above right] at (-2, 0) {$p_{4}$};
\draw[line width = 2pt, red] (0, 0) -- (-2, 0) {};

\draw (-2, 4.5) -- (0, 4.5) {};
\draw (-2, 4.2) -- (-2, 4.7) {};
\draw (0, 4.2) -- (0, 4.7) {};
\node[above] at (-1, 4.5) {$\Delta^{k} = \delta^{k} = 1/2$};

\end{tikzpicture}